\documentclass[10pt]{article}
\usepackage[english]{babel}
\usepackage{amsfonts}
\usepackage[dvips]{graphicx}

\begin{document}

\title{\bf A result of existence and uniqueness for a cavity driven flow. Analytical expression of the solution.}

\date{2009, December}

\author{\bf Gianluca Argentini \\
\normalsize{Research \& Development Dept., Riello Burners - Italy}\\
\normalsize gianluca.argentini@rielloburners.com \\
\normalsize gianluca.argentini@gmail.com \\}

\maketitle

\noindent{\bf Abstract}\\
In this work a result of existence and uniqueness for a plane cavity driven steady flow is deduced using an analytical method for the resolution of a linear partial differential problem on a triangular domain. The solution admits a symbolic expression based on integration over the domain. Some examples of flow are computed and graphed. In particular, it is shown a realistic example of a shear-forced flow with two eddies, usually computed only by numerical methods. The mathematical techniques used for the demonstration of the main result are elementary.\\


\newtheorem{theo}{Theorem}

\section{A theorem of existence and uniqueness}

Let $\Omega$ be (the interior of) a triangular domain in $\mathbb{R}^2$, the cartesian $\{x,y\}$ plane, with vertices $O=(0,0)$, $A=(2a,0)$, $B=(a,a)$, where $a$ is a positive real number. Note that the triangle $OAB$ is rectangular and $\overline{OB}=\overline{BA}$. Let be $f(x,y) \in C^0(\overline{\Omega},\mathbb{R})$. We want to resolve the differential problem

\begin{equation}\label{diffProblem}
	- \partial^2_{xx}\phi + \partial^2_{yy}\phi = f	\hspace{0.3cm} \textnormal{in} \hspace{0.1cm} \Omega, \hspace{0.3cm} \phi = 0 \hspace{0.3cm} \textnormal{on} \hspace{0.1cm} \partial\Omega
\end{equation}

\noindent for a function $\phi(x,y)$, $\phi \in C^2(\overline{\Omega},\mathbb{R})$. Note that the partial differential equation $- \partial^2_{xx}\phi + \partial^2_{yy}\phi = f$ admits a general solution of the form (see \cite{jeffrey} or \cite{sneddon}) $\phi(x,y)=g(-x+y)+h(x+y)+\phi_0(x,y)$, where $g$ and $h$ are arbitrary real functions and $\phi_0$ is a particular solution of the equation. But this general expression is not very useful for applying the boundary condition $\phi = 0$ on $\partial\Omega$. It is more interesting and instructive the following direct method.\\

\begin{figure}[ht!]
	\begin{center}
	\includegraphics[width=7cm]{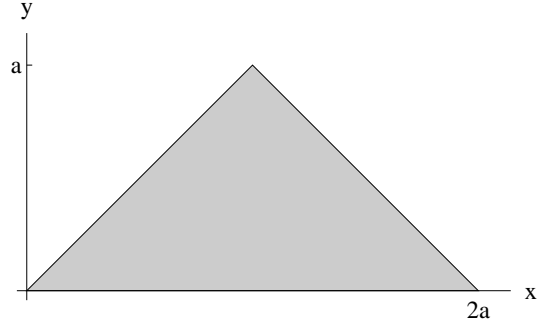}
	\caption{{\it The domain in the $xy$-plane}.}
	\label{trianglexy}
	\end{center}
\end{figure}

Consider the differential operator $- \partial^2_{xx} + \partial^2_{yy}$ written as $(\partial_x+\partial_y)(-\partial_x+\partial_y)$, and consider a linear transformation rule for cartesian coordinates $X=ax+by$, $Y=cx+dy$. If we want to have $2\partial_X = \partial_x+\partial_y$ and $2\partial_Y = -\partial_x+\partial_y$, using the chain rule it must be $a=b=1=d=1$ and $c=-1$, that is

\begin{equation}\label{transformationRule}
 X=x+y, \hspace{0.2cm} Y=-x+y
\end{equation}

\noindent The transformation is invertible:

\begin{equation}\label{inverseTransformationRule}
 2x=X-Y, \hspace{0.2cm} 2y=X+Y
\end{equation}

\noindent With the notation $\Phi(X,Y)=\phi(x(X,Y),y(X,Y))$ and analogous for $f$, the differential equation $- \partial^2_{xx}\phi + \partial^2_{yy}\phi = f$ becomes

\begin{equation}\label{diffEqnXY}
	4\partial^2_{XY}\Phi(X,Y) = F(X,Y)
\end{equation}

\noindent Note that the transformation (\ref{transformationRule}) is a $45^{\circ}$-rotation and a $\sqrt{2}$-dilation of the plane $\{x,y\}$. Also, the boundary condition doesn't change: $\Phi=0$ on $\partial\Omega$ (for simplicity we use for the domain in the plane $\{X,Y\}$ the same symbol $\Omega$ used for the plane $\{x,y\}$). For example, $\Phi(X,-X)=\phi(x,0)=0$. Therefore, the differential problem (\ref{diffProblem}) becomes

\begin{equation}\label{diffProblemXY}
	4\partial^2_{XY}\Phi = F	\hspace{0.3cm} \textnormal{in} \hspace{0.1cm} \Omega, \hspace{0.3cm} \Phi = 0 \hspace{0.3cm} \textnormal{on} \hspace{0.1cm} \partial\Omega
\end{equation}

\begin{figure}[ht!]
	\begin{center}
	\includegraphics[width=7cm]{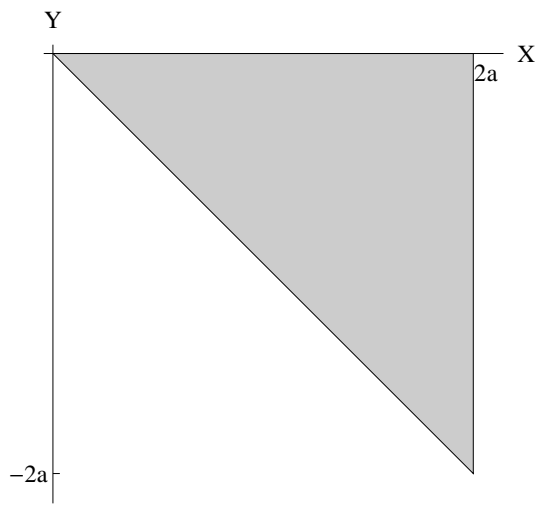}
	\caption{{\it The domain in the $XY$-plane}.}
	\label{triangleXYnew}
	\end{center}
\end{figure}

\noindent We remark the fact that the operator $4\partial^2_{XY}$ is the \textit{canonical form} of the differential operator $- \partial^2_{xx}\phi + \partial^2_{yy}$, which has the lines $y=x$ and $y=-x$ as \textit{characteristic curves} (\cite{jeffrey} or \cite{sneddon}).\\

Now we want to discuss the resolution of the differential problem (\ref{diffProblemXY}).

\begin{figure}[ht!]
	\begin{center}
	\includegraphics[width=6cm]{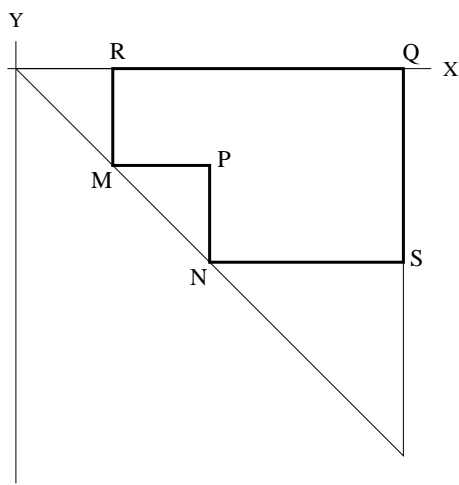}
	\caption{{\it Domain of integration}.}
	\label{domainIntXY}
	\end{center}
\end{figure}

\noindent Let $P=(X,Y)$ be a point in the interior of the domain $\Omega$. Then we can construct the polygon $\Sigma$ using segments parallel to $X$ and $Y$ axes (see Fig.\ref{domainIntXY}). Note that $M=(-Y,Y)$, $R=(-Y,0)$, $Q=(2a,0)$, $S=(2a,-X)$, $N=(X,-X)$. Using the identity 2$\partial^2_{XY}=(\partial_X\partial_Y + \partial_Y\partial_X)$, from the differential equation it follows that

\begin{equation}\label{eqnDiffIntegrated}
	2\int_{\Sigma} \left[ \partial_X\partial_Y\Phi(X,Y) + \partial_Y\partial_X\Phi(X,Y) \right] dXdY = \int_{\Sigma} F(X,Y) dXdY
\end{equation}

\noindent Now apply the Green theorem (\cite{matthews} or \cite{stewart}) to the first integral:

\begin{equation}
	\int_{\Sigma} \left(\partial_X\partial_Y\Phi + \partial_Y\partial_X\Phi\right) dXdY = \int_{\partial\Sigma}\left(\partial_Y\Phi \hspace{0.1cm} dY - \partial_X\Phi \hspace{0.1cm} dX\right)
\end{equation}	
	
\noindent It is now simple to calculate the line integral along the edges of the polygon $\Sigma$ (note that the boundary must be walked in counter-clockwise sense):

\begin{eqnarray}\label{lineIntegral}
	\int_{\partial\Sigma}\left(\partial_Y\Phi \hspace{0.1cm} dY - \partial_X\Phi \hspace{0.1cm} dX\right) = \nonumber \\
	= -2\Phi(P)+2\Phi(N)-2\Phi(S)+2\Phi(Q)-2\Phi(R)+2\Phi(M)
\end{eqnarray}

\noindent So we have

\begin{eqnarray}
	-2\Phi(P)+2\Phi(N)-2\Phi(S)+2\Phi(Q)-2\Phi(R)+2\Phi(M) = \nonumber \\
	= \int_{\Sigma} F dXdY
\end{eqnarray}

\noindent Now apply the boundary condition $\Phi_{|\partial\Omega} = 0$: it follows that

\begin{equation}\label{solution1}
	\Phi(P) = \Phi(X,Y) = -\frac{1}{2}\int_{\Sigma(X,Y)} F(t,s) \hspace{0.1cm} dt  \hspace{0.1cm} ds
\end{equation}

\noindent with the consequence that, if the point $P(X,Y)$ lies on the boundary of $\Omega$, that is if $P=M=N$ or $P=S$ or $P=R$, the function $F$ must satisfy the necessary condition

\begin{equation}
	0 = \int_{\Sigma(X,Y)} F(t,s) \hspace{0.1cm} dt  \hspace{0.1cm} ds \hspace{0.3cm} \forall (X,Y)\in\partial\Omega
\end{equation}

\noindent It is easy to see that previous condition can be written in a more explicit fashion:

\begin{equation}\label{necConditionF}
	\int_{X}^{2a}\int_{-X}^{0} F(t,s) \hspace{0.1cm} ds \hspace{0.1cm} dt = 0 \hspace{0.3cm} \forall X\in[0,2a]
\end{equation}

\noindent Therefore we have shown that a solution to differential problem (\ref{diffProblemXY}), and hence to (\ref{diffProblem}), exists if and only if $F$ satisfies condition (\ref{necConditionF}). Also, formula (\ref{solution1}) is an analytical expression for a solution. Note that, denoted by $T$ the point $(X,0)$, the integral can be divided into the two integrals defined on the two simple rectangles $PTRM$ and $NSQT$.\\

Now we discuss uniqueness of solution. Suppose to have two solutions $\Phi_1$ and $\Phi_2$ for the problem (\ref{diffProblemXY}). Then $\Phi=\Phi_1-\Phi_2$ is a function such that $\partial^2_{XY}\Phi=0$ $\forall (X,Y)\in\Omega$ and $\Phi_{|\partial\Omega}=0$. Note that we can write

\begin{equation}
	\partial_Y\left[\partial_X\Phi\right]^2 = 2  \hspace{0.1cm} \partial_X\Phi  \hspace{0.1cm} \partial_{XY}^2\Phi=0
\end{equation}

\noindent Applying the Green theorem to domain $\Sigma$ for the expression $\partial_Y\left[\partial_X\Phi\right]^2$, we have

\begin{eqnarray}\label{boundaryIntegrals}
	0 = \int_{\partial\Sigma}\left[\partial_X\Phi\right]^2dX = \nonumber \\
	= \int_{N}^{S}\left[\partial_X\Phi\right]^2dX + \int_{Q}^{R}\left[\partial_X\Phi\right]^2dX + \int_{M}^{P}\left[\partial_X\Phi\right]^2dX
\end{eqnarray}

\noindent Using integration by parts, the following identity holds:

\begin{equation}
	\int\left[\partial_X\Phi\right]^2dX = \Phi \hspace{0.1cm} \partial_X\Phi - \int\left[\Phi \hspace{0.1cm} \partial_{XX}^2\Phi\right]dX
\end{equation}

\noindent Hence, being $\Phi_{|\partial\Omega}=0$, the second integral in (\ref{boundaryIntegrals}) is null, therefore

\begin{equation}
	\int_{N}^{S}\left[\partial_X\Phi\right]^2dX + \int_{M}^{P}\left[\partial_X\Phi\right]^2dX = 0
\end{equation}

\noindent The two integrals are evaluated in the same sense of the integration path, so that $\partial_X\Phi=0$ along the segments $NS$ and $MP$, therefore $\Phi(P)=\Phi(M)$. But $\Phi(M)=0$, being $M\in\partial\Omega$, so for a generic point $P=(X,Y)$ we have $\Phi(P)=0$. The solutions $\Phi_1$ and $\Phi_2$ are identical.\\

We have shown the following result (remember that $X=x+y$, $Y=-x+y$):

\begin{theo}\label{theoremExistUniq}
	Let $f$ be a real function of $C^0(\overline{\Omega},\mathbb{R})$ such that
	\begin{equation}\label{theoCondition}
		\int_{X}^{2a}\int_{-X}^{0} f\left(\frac{t-s}{2},\frac{t+s}{2}\right) \hspace{0.1cm} ds \hspace{0.1cm} dt = 0 \hspace{0.3cm} \forall X\in[0,2a]
	\end{equation}
	Then the differential problem 
	\begin{equation}
		- \partial^2_{xx}\phi + \partial^2_{yy}\phi = f	\hspace{0.3cm} in \hspace{0.1cm} \Omega, \hspace{0.3cm} \phi = 0 \hspace{0.3cm} on \hspace{0.1cm} \partial\Omega
	\end{equation}
	has one and only one solution in the space $C^2(\overline{\Omega},\mathbb{R})$. The solution is given by the formula
	\begin{eqnarray}\label{formulaSolution}
		\phi(x,y) = &-&\frac{1}{2}\int_{x-y}^{x+y}\int_{-x+y}^0 f\left(\frac{t-s}{2},\frac{t+s}{2}\right) \hspace{0.1cm}ds \hspace{0.1cm} dt - \nonumber \\ 	
		&-&\frac{1}{2}\int_{x+y}^{2a}\int_{-x-y}^0 f\left(\frac{t-s}{2},\frac{t+s}{2}\right) \hspace{0.1cm}ds \hspace{0.1cm} dt
	\end{eqnarray}
\end{theo}

\section{An application: cavity driven flows}

In this section we discuss an application of previous theorem to a problem of two-dimensional cavity driven flow, that is a plane flow confined in a cavity and induced by the stress due to a primary flow external to the cavity (see \cite{shankar}). This phenomenon has great importance in scientific research (see e.g. \cite{madaniCom}) and technological applications. Assume that the cavity has the shape of the triangle $OAB$ of Fig.\ref{trianglexy} in the $xy$-plane. Stress due to the primary flow acts on the horizontal edge $OA$. We suppose that the fluid is newtonian and incompressible, that is plane stress $\mathbb{T}$ and plane strain-rate $\mathbb{D}$ tensors are linked by the formula (see \cite{gurtin})

\begin{equation}
	\mathbb{T} = 2\mu \mathbb{D}
\end{equation}

\noindent where $\mu$ is the dynamic viscosity and $2\mathbb{D}_{ij}=\left(\partial_j v_i + \partial_i v_j \right)$ (see \cite{madani}), where $\bf{v}$=$(v_1,v_2)$=$(v_x,v_y)$ is the flow velocity field. Plane incompressible flows admit a stream function (\cite{madani}), that is a function $\Psi(x,y)$ such that

\begin{equation}
	u = v_x = \partial_y \Psi, \hspace{0.3cm} v = v_y = - \partial_x \Psi
\end{equation}

\noindent Therefore, a plane newtonian incompressible flow is described by the partial differential equation

\begin{equation}\label{partialDiffPsi}
	- \partial_{xx}^2 \Psi + \partial_{yy} \Psi = \frac{1}{\mu}\mathbb{T}_{xy}
\end{equation}

In the next of the paper we suppose to know the analytical expression of $\mathbb{T}_{xy}$ and we try to find a solution of (\ref{partialDiffPsi}) for a stream function $\Psi$ such that $\Psi_{|\partial \Omega} = 0$. This boundary condition is usual for plane incompressible flow (see \cite{madani} and \cite{erturk}), but in the case of a cavity driven flow it (or an analogous $\Psi_{|\partial \Omega}$ = const) has an important physical meaning. In fact, if $\Psi_{|\partial \Omega} = 0$, then $\partial\Omega$ is a level curve for $\Psi$, therefore at each point of the boundary $\nabla\Psi$ is orthogonal to the tangent of the boundary itself (\cite{stewart}). But $\nabla\Psi = (-v,u)$, which is orthogonal to the flow velocity field $(u,v)$. Therefore, at each point of $\partial\Omega$, the geometrical tangent and the velocity field are parallel, that is the flow is confined into the cavity $\Omega$.\\
\noindent Applying theorem (\ref{theoremExistUniq}), it can be stated that if $\mathbb{T}_{xy} \in C^0(\overline{\Omega},\mathbb{R})$, then there is a unique stream function $\Psi \in C^2(\overline{\Omega},\mathbb{R})$ solving the linear equation (\ref{partialDiffPsi}) with boundary condition $\Psi|_{\partial\Omega}=0$. Note that \cite{albensoeder} consider a nonlinear problem about cavity driven flow where uniqueness can fail.

We consider at first the more simple analytical form for a possible stress:

\begin{equation}
	\mathbb{T}_{xy} = \mu(c_1y+c_2)
\end{equation}

\noindent In this case, along the horizontal edge $OA$ ($y=0$) of the cavity the stress is constant. From theorem (\ref{theoremExistUniq}), a solution to our differential problem exists if the function $c_1y+c_2$ satisfies the condition (\ref{theoCondition}). It is easy to show that the condition is satisfied for all $X\in[0,2a]$ if and only if $2 c_2 = - a c_1$. Note that in this case the stress has expression

\begin{equation}
	\mathbb{T}_{xy} = \mu c_2 \left(-\frac{2}{a}y+1\right)
\end{equation}

\noindent and for $y=\frac{a}{2}$ it changes its sign. So flow can recirculate. We make the choice $c_2=-8a$, so that $c_1=16$. The differential problem to solve is

\begin{equation}\label{diffProblemTest}
	- \partial^2_{xx}\phi + \partial^2_{yy}\phi = 16y-8a	\hspace{0.3cm} \textnormal{in} \hspace{0.1cm} \Omega, \hspace{0.3cm} \phi = 0 \hspace{0.3cm} \textnormal{on} \hspace{0.1cm} \partial\Omega
\end{equation}

\begin{figure}[ht!]
	\begin{center}
	\includegraphics[width=6cm]{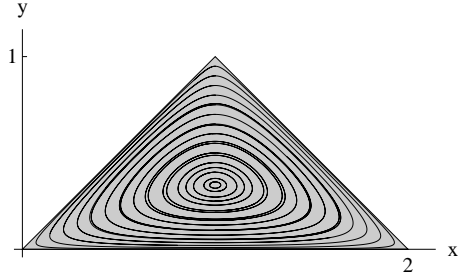}
	\caption{{\it Flow path-lines in the case a}=1.}
	\label{pathLines}
	\end{center}
\end{figure}

\noindent From formula (\ref{formulaSolution}), using the transformation rule $X=x+y$, $Y=-x+y$ the solution to (\ref{diffProblemTest}) can be computed by

\begin{eqnarray}\label{formulaSolutionTest}
		\Psi(x,y) = &-&4\int_{x-y}^{x+y}\int_{-x+y}^0 (t+s-a) \hspace{0.1cm}ds \hspace{0.1cm} dt - \nonumber \\ 	
		&-&4\int_{x+y}^{2a}\int_{-x-y}^0 (t+s-a) \hspace{0.1cm}ds \hspace{0.1cm} dt
\end{eqnarray}

\noindent which gives the expression

\begin{equation}
	\Psi(x,y)=2y^3-2x^2y-4ay^2+4axy
\end{equation}

\noindent for the stream function of the flow. The velocity field is $(\partial_y\Psi,-\partial_x\Psi)=(-2x^2+6y^2+4ax-8ay,4xy-4ay)$. It is interesting to find the points where the velocity is null. Solving the algebraic system $(-2x^2+6y^2+4ax-8ay,4xy-4ay)=(0,0)$, we find as expected the three vertices $(0,0)$, $(2a,0)$ and $(a,a)$, and also the interior point $(\frac{a}{3},a)$ which is the center of the recirculation gyre (see Fig.\ref{pathLines}).

\begin{figure}[ht!]
	\begin{center}
	\includegraphics[width=6cm]{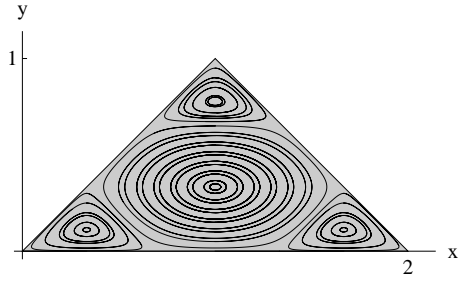}
	\caption{{\it Flow path-lines in the case of the sinusoidal stress, with a}=1, A=5.}
	\label{sinusoidalPathLines}
	\end{center}
\end{figure}

\begin{figure}[ht!]
	\begin{center}
	\includegraphics[width=6cm]{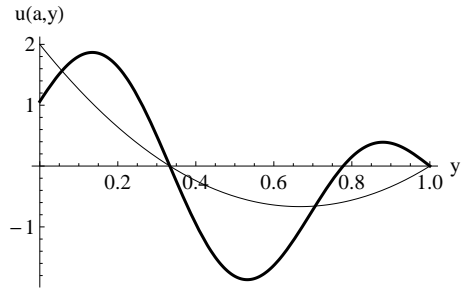}
	\caption{{\it Comparison of $u=\partial_y\Psi$ in the linear (thin line) and sinusoidal (thick line) case}.}
	\label{graphsSin}
	\end{center}
\end{figure}

Now we consider a more interesting case. Let the stress be described by a sinusoidal expression of the form 

\begin{equation}
	\mathbb{T}_{xy} = A \mu \hspace{0.1cm} \textnormal{cos}(ky)
\end{equation}

\noindent with $A$ and $k$ real numbers. Using (\ref{theoCondition}), it is easy to show that if we suppose $q=0$, then

\begin{equation}
	k = m\frac{\pi}{a},	\hspace{0.3cm} m=2n+1, \hspace{0.3cm} n \in \mathbb{N}
\end{equation}

\noindent is the condition for existence and uniqueness of a flow in the triangular cavity. Consider $m=1$. By integration (\ref{formulaSolution}), the analytic form of the stream function, solution of the differential problem, is

\begin{equation}
	\Psi=-\frac{2Aa^2}{9\pi^2}\left[\textnormal{cos}\left(\frac{3\pi}{a}y\right)+\textnormal{cos}\left(\frac{3\pi}{2a}(x-y)\right)-2\textnormal{cos}^2\left(\frac{3\pi}{4a}(x+y)\right)\right]
\end{equation}

\noindent and Fig.\ref{sinusoidalPathLines} shows some path lines, where one primary central eddy and three secondary eddies are present. It is also interesting to draw the graph of $u=\partial_y\Psi$ for $x=a$ and $y$ variable in the range $[0,a]$: there are two values of $y$, not equal to $a$, for which $u=0$. One of the two values is equal to the value where $u=0$ in the previous case of the linear stress (see Fig.\ref{graphsSin}).\\

A more realistic example is based on a stress with analytical expression of the form

\begin{equation}
	\mathbb{T}_{xy}=\mu \sum_{m,n=0}^4a_{m,n} \hspace{0.1cm} x^my^n
\end{equation}

\noindent Applying condition (\ref{theoCondition}) for the computation of the coefficients $a_{m,n}$, a possible stream function is

\begin{equation}
	\Psi(x,y)=(2y^3-2x^2y-4ay^2+4axy)(y-100x^2-a)\left(y+\frac{1}{4}x-\frac{5}{6}a\right)
\end{equation}

\begin{figure}[ht!]
	\begin{center}
	\includegraphics[width=6cm]{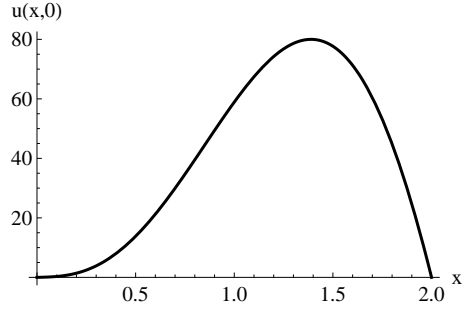}
	\caption{{\it Profile of $u(x,0)=\partial_y\Psi(x,0)$ in the case of a variable shear stress}.}
	\label{ugraphRealistic}
	\end{center}
\end{figure}

\noindent The horizontal component $u$ of the velocity field, along the $x$-axes, is a 5-order $x$-polynomial whose graph is shown on Fig.\ref{ugraphRealistic}. The stress acts on the horizontal segment of the triangular domain as a variable shear of positive sign.

\begin{figure}[ht!]
	\begin{center}
	\includegraphics[width=6cm]{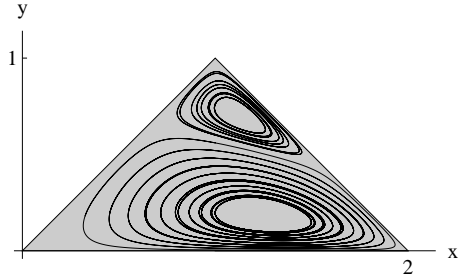}
	\caption{{\it Flow path-lines in the case of a variable shear stress, with a}=1.}
	\label{realisticPathLines}
	\end{center}
\end{figure}

\noindent The resulting flow path-lines (see Fig.\ref{realisticPathLines}) show the presence of a primary gyre, and of a secondary gyre near the vertex opposite to the edge subjected to external stress. This image is similar to a corrisponding picture (fig.2(b)) in \cite{erturk}, where stream function is computed by numerical method.

\end{document}